\documentclass[a4paper,11pt]{amsart}
\usepackage{amsmath,amssymb,amsfonts,latexsym}

\newtheorem{theorem}{Theorem}[section]

\newtheorem{corollary}[theorem]{Corollary}

\input{tcilatex}

\begin{document}
\title[\textbf{On the Dirichlet's type of Eulerian polynomials}]{\textbf{On
the Dirichlet's type of Eulerian polynomials}}
\author[\textbf{S. Araci}]{\textbf{Serkan Araci}}
\address{\textbf{University of Gaziantep, Faculty of Science and Arts,
Department of Mathematics, 27310 Gaziantep, TURKEY}}
\email{\textbf{mtsrkn@hotmail.com}}
\author[\textbf{M. Acikgoz}]{\textbf{Mehmet Acikgoz}}
\address{\textbf{University of Gaziantep, Faculty of Science and Arts,
Department of Mathematics, 27310 Gaziantep, TURKEY}}
\email{\textbf{acikgoz@gantep.edu.tr}}
\author[\textbf{D. Gao}]{\textbf{Deyao Gao}}
\address{\textbf{Yuren Lab., No. 8 Tongsheng Road, Changsha, P. R. China}}
\email{\textbf{13607433711@163.com}}

\begin{abstract}
In the present paper, we introduce Eulerian polynomials attached to $\chi $
by using $p$-adic $q$-integral on $%
%TCIMACRO{\U{2124} }%
%BeginExpansion
\mathbb{Z}
%EndExpansion
_{p}$. Also, we give new interesting identities via the generating functions
of Dirichlet's type of Eulerian polynomials. After, by applying Mellin
transformation to this generating function of Dirichlet's type of Eulerian
polynomials, we derive $L$-function for Eulerian polynomials which
interpolates of Dirichlet's type of Eulerian polynomials at negative
integers.

\vspace{2mm}\noindent \textsc{2010 Mathematics Subject Classification.}
11S80, 11B68.

\vspace{2mm}

\end{abstract}

\maketitle

%%%%%%%%%%%%%%%%%%%%%%%%%%%%%%%%%%%%%%%%%%%%%%%%%%%%%%%%%%%%%%%%%%%

%%%%%%%%%%%%%%%%%%%%%%%%%%%%%%%%%%%%%%%%%%%%%%%%%%%%%%%%%%%%%%%%%%%

%%%%%%%%%%%%%%%%%%%%%%%%%%%%%%%%%%%%%%%%%%%%%%%%%%%%%%%%%%%%%%%%%%%

\section{\textbf{Introduction}}

%%%%%%%%%%%%%%%%%%%%%%%%%%%%%%%%%%%%%%%%%%%%%%%%%%%%%%%%%%%%%%%%%%%

Recently, Kim $et$ $al.$ have studied to Eulerian polynomials. They gave not
only Witt's formula for Eulerian polynomials but also relations between
Genocchi, Tangent and Euler numbers (for more details, see \cite{KIM}). In
arithmetic works of T. Kim have introduced many different generating
functions of families of Bernoulli, Euler, Genocchi numbers and polynomials
by using $p$-adic $q$-integral on $%
%TCIMACRO{\U{2124} }%
%BeginExpansion
\mathbb{Z}
%EndExpansion
_{p}$ (see [1-5]). After, many mathematicians are motivated from his papers
and introduced new generating function for special functions (for more
information about this subject, see [16-28]). Y. Simsek also gave new $q$%
-twisted Euler numbers and polynomials and ($h,q$)-Bernoulli numbers and
polynomials by using Kim's $p$-adic $q$-integral on $%
%TCIMACRO{\U{2124} }%
%BeginExpansion
\mathbb{Z}
%EndExpansion
_{p}$. He also derived some interesting properties in his works \cite{Simsek}%
, \cite{Simsek 1}, \cite{Simsek 2}.

The $p$-adic $q$-integral on $%
%TCIMACRO{\U{2124} }%
%BeginExpansion
\mathbb{Z}
%EndExpansion
_{p}$ was originally defined by Kim. He also investigated that $p$-adic $q$%
-integral on $%
%TCIMACRO{\U{2124} }%
%BeginExpansion
\mathbb{Z}
%EndExpansion
_{p}$ is related to non-Archimedean combinatorial analysis in mathematical
physics. That is, the functional equation of the $q$-zeta function, the $q$%
-Stirling numbers and $q$-Mahler theory and so on (for details, see\cite%
{KIM4}, \cite{KIM5}).

We firstly list some properties of familiar Eulerian polynomials for sequel
of this paper as follows:

As it is well-known, the Eulerian polynomials, $\mathcal{A}_{n}\left(
x\right) $ are given by means of the following generating function:%
\begin{equation}
e^{\mathcal{A}\left( x\right) t}=\sum_{n=0}^{\infty }\mathcal{A}_{n}\left(
x\right) \frac{t^{n}}{n!}=\frac{1-x}{e^{t\left( 1-x\right) }-x}
\label{equation 19}
\end{equation}

where $\mathcal{A}^{n}\left( x\right) :=\mathcal{A}_{n}\left( x\right) $ as
symbolic. To find Eulerian polynomials, it has the following recurrence
relation:%
\begin{equation}
\left( \mathcal{A}\left( t\right) +\left( t-1\right) \right) ^{n}-t\mathcal{A%
}_{n}\left( t\right) =\left\{ 
\begin{array}{cc}
1-t & \text{if }n=0 \\ 
0 & \text{if }n\neq 0,%
\end{array}%
\right.  \label{equation 20}
\end{equation}

(for details, see \cite{KIM}).

Suppose that $p$ be a fixed odd prime number. Throughout this paper, we use
the following notations. By $%
%TCIMACRO{\U{2124} }%
%BeginExpansion
\mathbb{Z}
%EndExpansion
_{p}$, we denote the ring of $p$-adic rational integers, $%
%TCIMACRO{\U{211a} }%
%BeginExpansion
\mathbb{Q}
%EndExpansion
$ denotes the field of rational numbers, $%
%TCIMACRO{\U{211a} }%
%BeginExpansion
\mathbb{Q}
%EndExpansion
_{p}$ denotes the field of $p$-adic rational numbers, and $%
%TCIMACRO{\U{2102} }%
%BeginExpansion
\mathbb{C}
%EndExpansion
_{p}$ denotes the completion of algebraic closure of $%
%TCIMACRO{\U{211a} }%
%BeginExpansion
\mathbb{Q}
%EndExpansion
_{p}$. Let $%
%TCIMACRO{\U{2115} }%
%BeginExpansion
\mathbb{N}
%EndExpansion
$ be the set of natural numbers and $%
%TCIMACRO{\U{2115} }%
%BeginExpansion
\mathbb{N}
%EndExpansion
^{\ast }=%
%TCIMACRO{\U{2115} }%
%BeginExpansion
\mathbb{N}
%EndExpansion
\cup \left\{ 0\right\} $.

The $p$-adic absolute value is defined by 
\begin{equation*}
\left\vert p\right\vert _{p}=\frac{1}{p}\text{.}
\end{equation*}%
In this paper we assume $\left\vert q-1\right\vert _{p}<1$ as an
indeterminate. Let $UD\left( 
%TCIMACRO{\U{2124} }%
%BeginExpansion
\mathbb{Z}
%EndExpansion
_{p}\right) $ be the space of uniformly differentiable functions on $%
%TCIMACRO{\U{2124} }%
%BeginExpansion
\mathbb{Z}
%EndExpansion
_{p}$. For a positive integer $d$ with $\left( d,p\right) =1$, set 
\begin{eqnarray*}
X &=&X_{m}=\lim_{\overleftarrow{m}}%
%TCIMACRO{\U{2124} }%
%BeginExpansion
\mathbb{Z}
%EndExpansion
/dp^{m}%
%TCIMACRO{\U{2124} }%
%BeginExpansion
\mathbb{Z}
%EndExpansion
\text{,} \\
X^{\ast } &=&\underset{\underset{\left( a,p\right) =1}{0<a<dp}}{\cup }a+dp%
%TCIMACRO{\U{2124} }%
%BeginExpansion
\mathbb{Z}
%EndExpansion
_{p}
\end{eqnarray*}

and%
\begin{equation*}
a+dp^{m}%
%TCIMACRO{\U{2124} }%
%BeginExpansion
\mathbb{Z}
%EndExpansion
_{p}=\left\{ x\in X\mid x\equiv a\left( \func{mod}dp^{m}\right) \right\} 
\text{,}
\end{equation*}

where $a\in 
%TCIMACRO{\U{2124} }%
%BeginExpansion
\mathbb{Z}
%EndExpansion
$ satisfies the condition $0\leq a<dp^{m}$.

Firstly, for introducing fermionic $p$-adic $q$-integral, we need some basic
information which we state here. A measure on $%
%TCIMACRO{\U{2124} }%
%BeginExpansion
\mathbb{Z}
%EndExpansion
_{p}$ with values in a $p$-adic Banach space $B$ is a continuous linear map%
\newline
\begin{equation*}
f\mapsto \int f(x)\mu =\int_{%
%TCIMACRO{\U{2124} }%
%BeginExpansion
\mathbb{Z}
%EndExpansion
_{p}}f(x)\mu (x)
\end{equation*}

from $C^{0}(%
%TCIMACRO{\U{2124} }%
%BeginExpansion
\mathbb{Z}
%EndExpansion
_{p},%
%TCIMACRO{\U{2102} }%
%BeginExpansion
\mathbb{C}
%EndExpansion
_{p})$, (continuous function on $%
%TCIMACRO{\U{2124} }%
%BeginExpansion
\mathbb{Z}
%EndExpansion
_{p}$ ) to $B$. We know that the set of locally constant functions from $%
%TCIMACRO{\U{2124} }%
%BeginExpansion
\mathbb{Z}
%EndExpansion
_{p}$ to $%
%TCIMACRO{\U{211a} }%
%BeginExpansion
\mathbb{Q}
%EndExpansion
_{p}$ is dense in $C^{0}(%
%TCIMACRO{\U{2124} }%
%BeginExpansion
\mathbb{Z}
%EndExpansion
_{p},%
%TCIMACRO{\U{2102} }%
%BeginExpansion
\mathbb{C}
%EndExpansion
_{p})$ so.\newline

Explicitly, for all $f\in C^{0}(%
%TCIMACRO{\U{2124} }%
%BeginExpansion
\mathbb{Z}
%EndExpansion
_{p},%
%TCIMACRO{\U{2102} }%
%BeginExpansion
\mathbb{C}
%EndExpansion
_{p})$, the locally constant functions\newline
\begin{equation*}
f_{n}=\sum_{i=0}^{p^{m}-1}f(i)1_{i+p^{m}%
%TCIMACRO{\U{2124} }%
%BeginExpansion
\mathbb{Z}
%EndExpansion
_{p}}\rightarrow f\text{ in }C^{0}
\end{equation*}

Now, set $\mu (i+p^{m}%
%TCIMACRO{\U{2124} }%
%BeginExpansion
\mathbb{Z}
%EndExpansion
_{p})=\int_{%
%TCIMACRO{\U{2124} }%
%BeginExpansion
\mathbb{Z}
%EndExpansion
_{p}}1_{i+p^{m}%
%TCIMACRO{\U{2124} }%
%BeginExpansion
\mathbb{Z}
%EndExpansion
_{p}}\mu $. Then $\int_{%
%TCIMACRO{\U{2124} }%
%BeginExpansion
\mathbb{Z}
%EndExpansion
_{p}}f\mu $, is given by the following Riemannian sum%
\begin{equation*}
\int_{%
%TCIMACRO{\U{2124} }%
%BeginExpansion
\mathbb{Z}
%EndExpansion
_{p}}f\mu =\lim_{m\rightarrow \infty }\sum_{i=0}^{p^{m}-1}f(i)\mu {(i+p^{m}%
%TCIMACRO{\U{2124} }%
%BeginExpansion
\mathbb{Z}
%EndExpansion
_{p})}
\end{equation*}

The following $q$-Haar measure is defined by Kim in \cite{KIM2} and \cite%
{KIM4}:%
\begin{equation*}
\mu _{q}(a+p^{m}%
%TCIMACRO{\U{2124} }%
%BeginExpansion
\mathbb{Z}
%EndExpansion
_{p})=\frac{q^{a}}{[p^{m}]_{q}}
\end{equation*}

So, for $f\in UD\left( 
%TCIMACRO{\U{2124} }%
%BeginExpansion
\mathbb{Z}
%EndExpansion
_{p}\right) $, the $p$-adic $q$-integral on $%
%TCIMACRO{\U{2124} }%
%BeginExpansion
\mathbb{Z}
%EndExpansion
_{p}$ is defined by Kim as follows:%
\begin{equation}
I_{q}\left( f\right) =\int_{%
%TCIMACRO{\U{2124} }%
%BeginExpansion
\mathbb{Z}
%EndExpansion
_{p}}f\left( \eta \right) d\mu _{q}\left( \eta \right) =\lim_{n\rightarrow
\infty }\frac{1}{\left[ p^{n}\right] _{q}}\sum_{\eta =0}^{p^{n}-1}q^{\eta
}f\left( \eta \right) \text{.}  \label{equation 2}
\end{equation}%
The bosonic integral is considered as the bosonic limit $q\rightarrow 1,$ $%
I_{1}\left( f\right) =\lim_{q\rightarrow 1}I_{q}\left( f\right) $. In \cite%
{Kim 2}, \cite{Kim 3} and \cite{Kim 4}, similarly, the $p$-adic fermionic
integration on $%
%TCIMACRO{\U{2124} }%
%BeginExpansion
\mathbb{Z}
%EndExpansion
_{p}$ \ defined by Kim as follows:%
\begin{equation}
I_{-q}\left( f\right) =\lim_{q\rightarrow -q}I_{q}\left( f\right) =\int_{%
%TCIMACRO{\U{2124} }%
%BeginExpansion
\mathbb{Z}
%EndExpansion
_{p}}f\left( x\right) d\mu _{-q}\left( x\right) \text{.}  \label{equation 3}
\end{equation}

By (\ref{equation 3}), we have the following well-known integral eguation:%
\begin{equation}
q^{n}I_{-q}\left( f_{n}\right) +\left( -1\right) ^{n-1}I_{-q}\left( f\right)
=\left[ 2\right] _{q}\sum_{l=0}^{n-1}\left( -1\right) ^{n-1-l}q^{l}f\left(
l\right)  \label{equation 4}
\end{equation}

Here $f_{n}\left( x\right) :=f\left( x+n\right) $. By (\ref{equation 4}), we
have the following equalities:

If $n$ odd, then%
\begin{equation}
q^{n}I_{-q}\left( f_{n}\right) +I_{-q}\left( f\right) =\left[ 2\right]
_{q}\sum_{l=0}^{n-1}\left( -1\right) ^{l}q^{l}f\left( l\right) \text{.}
\label{equation 5}
\end{equation}

If $n$ even, then we have%
\begin{equation}
I_{-q}\left( f\right) -q^{n}I_{-q}\left( f_{n}\right) =\left[ 2\right]
_{q}\sum_{l=0}^{n-1}\left( -1\right) ^{l}q^{l}f\left( l\right) \text{.}
\label{equation 6}
\end{equation}

Substituting $n=1$ into (\ref{equation 5}), we readily see the following%
\begin{equation}
qI_{-q}\left( f_{1}\right) +I_{-q}\left( f\right) =\left[ 2\right]
_{q}f\left( 0\right) \text{.}  \label{equation 7}
\end{equation}

Replacing $q$ by $q^{-1}$ in (\ref{equation 7}),\ we easily derive the
following%
\begin{equation}
I_{-q^{-1}}\left( f_{1}\right) +qI_{-q^{-1}}\left( f\right) =\left[ 2\right]
_{q}f\left( 0\right) \text{.}  \label{equation 8}
\end{equation}

In \cite{KIM}, Kim $et$ $al.$ is considered $f(x)=e^{-x\left( 1+q\right) t}$
in (\ref{equation 8}), then they gave Witt's formula of Eulerian polynomials
as follows:

For $n\in 
%TCIMACRO{\U{2115} }%
%BeginExpansion
\mathbb{N}
%EndExpansion
^{\ast }$,%
\begin{equation}
I_{-q^{-1}}\left( x^{n}\right) =\frac{\left( -1\right) ^{n}}{\left(
1+q\right) ^{n}}\mathcal{A}_{n}\left( -q\right) \text{.}  \label{equation 15}
\end{equation}

Now also, we consider $I_{-q^{-1}}\left( \chi \left( x\right) x^{n}\right) $
in the next section. We shall call as Dirichlet's type of Eulerian
polynomials. After we shall give arithmetic properties for Dirichlet's type
of Eulerian polynomials.

\section{\textbf{On the Dirichlet's type of Eulerian polynomials}}

Firstly, we consider the following equality by using (\ref{equation 5}): For 
$d$ odd natural numbers,%
\begin{gather}
\int_{%
%TCIMACRO{\U{2124} }%
%BeginExpansion
\mathbb{Z}
%EndExpansion
_{p}}f\left( x+d\right) d\mu _{-q^{-1}}\left( x\right) +q^{d}\int_{%
%TCIMACRO{\U{2124} }%
%BeginExpansion
\mathbb{Z}
%EndExpansion
_{p}}f\left( x\right) d\mu _{-q^{-1}}\left( x\right)  \label{equation 9} \\
=\left[ 2\right] _{q}\sum_{0\leq l\leq d-1}\left( -1\right)
^{l}q^{d-l+1}f\left( l\right) \text{.}  \notag
\end{gather}

Let $\chi $ be a Dirichlet's character of conductor $d,$ which is any
multiple of $p$ (=$odd$). Then, substituting $f(x)=\chi \left( x\right)
e^{-x\left( 1+q\right) t}$ in (\ref{equation 9}), then we compute as follows:%
\begin{gather*}
\int_{%
%TCIMACRO{\U{2124} }%
%BeginExpansion
\mathbb{Z}
%EndExpansion
_{p}}\chi \left( x+d\right) e^{-\left( x+d\right) \left( 1+q\right) t}d\mu
_{-q^{-1}}\left( x\right) +q^{d}\int_{%
%TCIMACRO{\U{2124} }%
%BeginExpansion
\mathbb{Z}
%EndExpansion
_{p}}\chi \left( x\right) e^{-x\left( 1+q\right) t}d\mu _{-q^{-1}}\left(
x\right) \\
=\left[ 2\right] _{q}\sum_{0\leq l\leq d-1}\left( -1\right)
^{l}q^{d-l+1}\chi \left( l\right) e^{-l\left( 1+q\right) t}
\end{gather*}

After some applications, we discover the following%
\begin{equation}
\int_{%
%TCIMACRO{\U{2124} }%
%BeginExpansion
\mathbb{Z}
%EndExpansion
_{p}}\chi \left( x\right) e^{-x\left( 1+q\right) t}d\mu _{-q^{-1}}\left(
x\right) =\left[ 2\right] _{q}\sum_{l=0}^{d-1}\left( -1\right)
^{l}q^{d-l+1}\chi \left( l\right) \frac{e^{-l\left( 1+q\right) t}}{%
e^{-d\left( 1+q\right) t}+q^{d}}\text{.}  \label{equation 10}
\end{equation}

Let $\mathcal{F}_{q}\left( t\mid \chi \right) =\sum_{n=0}^{\infty }\mathcal{A%
}_{n,\chi }\left( -q\right) \frac{t^{n}}{n!}$. Then, we introduce the
following definition of generating function of Dirichlet's type of Eulerian
polynomials.

\begin{definition}
For $n\in 
%TCIMACRO{\U{2115} }%
%BeginExpansion
\mathbb{N}
%EndExpansion
^{\ast }$, then we define the following:%
\begin{equation}
\sum_{n=0}^{\infty }\mathcal{A}_{n,\chi }\left( -q\right) \frac{t^{n}}{n!}=%
\left[ 2\right] _{q}\sum_{l=0}^{d-1}\left( -1\right) ^{l}q^{d-l+1}\chi
\left( l\right) \frac{e^{-l\left( 1+q\right) t}}{e^{-d\left( 1+q\right)
t}+q^{d}}\text{.}  \label{equation 11}
\end{equation}
\end{definition}

By (\ref{equation 10}) and (\ref{equation 11}), we state the following
theorem which is the Witt's formula for Dirichlet's type of Eulerian
polynomials.

\begin{theorem}
The following identity holds true:%
\begin{equation}
I_{-q^{-1}}\left( \chi \left( x\right) x^{n}\right) =\frac{\left( -1\right)
^{n}}{\left( 1+q\right) ^{n}}\mathcal{A}_{n,\chi }\left( -q\right) \text{.}
\label{equation 14}
\end{equation}
\end{theorem}

By using (\ref{equation 11}), we discover the following applications:%
\begin{eqnarray*}
\sum_{n=0}^{\infty }\mathcal{A}_{n,\chi }\left( -q\right) \frac{t^{n}}{n!}
&=&\left[ 2\right] _{q}\sum_{l=0}^{d-1}\left( -1\right) ^{l}q^{d-l+1}\chi
\left( l\right) \frac{e^{-l\left( 1+q\right) t}}{e^{-d\left( 1+q\right)
t}+q^{d}} \\
&=&\left[ 2\right] _{q}\sum_{l=0}^{d-1}\left( -1\right) ^{l}q^{-l+1}\chi
\left( l\right) e^{-l\left( 1+q\right) t}\sum_{m=0}^{\infty }\left(
-1\right) ^{m}q^{-md}e^{-md\left( 1+q\right) t} \\
&=&q\left[ 2\right] _{q}\sum_{m=0}^{\infty }\sum_{l=0}^{d-1}\left( -1\right)
^{l+md}\chi \left( l+md\right) q^{-\left( l+md\right) }e^{-\left(
l+md\right) \left( 1+q\right) t} \\
&=&q\left[ 2\right] _{q}\sum_{m=0}^{\infty }\left( -1\right) ^{m}\chi \left(
m\right) q^{-m}e^{-m\left( 1+q\right) t}\text{.}
\end{eqnarray*}

Thus, we get the following theorem.

\begin{theorem}
The following%
\begin{equation}
\mathcal{F}_{q}\left( t\mid \chi \right) =\sum_{n=0}^{\infty }\mathcal{A}%
_{n,\chi }\left( -q\right) \frac{t^{n}}{n!}=\left[ 2\right]
_{q}\sum_{m=0}^{\infty }\frac{\left( -1\right) ^{m}\chi \left( m\right)
q^{-m}e^{-m\left( 1+q\right) t}}{q^{m-1}}  \label{equation 12}
\end{equation}%
is true.
\end{theorem}

By considering Taylor expansion of $e^{-m\left( 1+q\right) t}$ in (\ref%
{equation 12}), we procure the following theorem.

\begin{theorem}
For $n\in 
%TCIMACRO{\U{2115} }%
%BeginExpansion
\mathbb{N}
%EndExpansion
$, then we have%
\begin{equation}
\frac{\left( -1\right) ^{n}}{q\left( 1+q\right) ^{n+1}}\mathcal{A}_{n,\chi
}\left( -q\right) =\sum_{m=1}^{\infty }\frac{\left( -1\right) ^{m}\chi
\left( m\right) m^{n}}{q^{m}}\text{.}  \label{equation 13}
\end{equation}
\end{theorem}

From (\ref{equation 14}) and (\ref{equation 13}), we easily derive the
following corollary:

\begin{corollary}
For $n\in 
%TCIMACRO{\U{2115} }%
%BeginExpansion
\mathbb{N}
%EndExpansion
$, then we procure the following%
\begin{equation*}
\lim_{m\rightarrow \infty }\sum_{x=1}^{p^{m}-1}\frac{\left( -1\right)
^{x}\chi \left( x\right) x^{n}}{q^{x}}=2\sum_{m=1}^{\infty }\frac{\left(
-1\right) ^{m}\chi \left( m\right) m^{n}}{q^{m-2}}\text{.}
\end{equation*}
\end{corollary}

Now, we give distribution formula for Dirichlet's type of Eulerian
polynomials by using $p$-adic $q$-integral on $%
%TCIMACRO{\U{2124} }%
%BeginExpansion
\mathbb{Z}
%EndExpansion
_{p}$, as follows:%
\begin{eqnarray*}
\int_{%
%TCIMACRO{\U{2124} }%
%BeginExpansion
\mathbb{Z}
%EndExpansion
_{p}}\chi \left( x\right) x^{n}d\mu _{-q^{-1}}\left( x\right)
&=&\lim_{m\rightarrow \infty }\frac{1}{\left[ dp^{m}\right] _{-q^{-1}}}%
\sum_{x=0}^{dp^{m}-1}\left( -1\right) ^{x}\chi \left( x\right) x^{n}q^{-x} \\
&=&\frac{d^{n}}{\left[ d\right] _{-q^{-1}}}\sum_{a=0}^{d-1}\left( -1\right)
^{a}\chi \left( a\right) q^{-a}\left( \lim_{m\rightarrow \infty }\frac{1}{%
\left[ p^{m}\right] _{-q^{-d}}}\sum_{x=0}^{p^{m}-1}\left( -1\right)
^{x}\left( \frac{a}{d}+x\right) ^{n}q^{-dx}\right) \\
&=&\frac{d^{n}}{\left[ d\right] _{-q^{-1}}}\sum_{a=0}^{d-1}\left( -1\right)
^{a}\chi \left( a\right) q^{-a}\int_{%
%TCIMACRO{\U{2124} }%
%BeginExpansion
\mathbb{Z}
%EndExpansion
_{p}}\left( \frac{a}{d}+x\right) ^{n}d\mu _{-q^{-d}}\left( x\right) .
\end{eqnarray*}

Thus, we state the following theorem.

\begin{theorem}
The following identity holds true:%
\begin{equation}
\frac{\left( -1\right) ^{n}}{\left( 1+q\right) ^{n}}\mathcal{A}_{n,\chi
}\left( -q\right) =\frac{d^{n}}{\left[ d\right] _{-q^{-1}}}%
\sum_{a=0}^{d-1}\left( -1\right) ^{a}\chi \left( a\right) q^{-a}\int_{%
%TCIMACRO{\U{2124} }%
%BeginExpansion
\mathbb{Z}
%EndExpansion
_{p}}\left( \frac{a}{d}+x\right) ^{n}d\mu _{-q^{-d}}\left( x\right) \text{.}
\label{equation 16}
\end{equation}
\end{theorem}

From this, we notice that the above equation is related to $q$-Genocchi
polynomials with weight zero, $\widetilde{G}_{n,q}\left( x\right) $, and $q$%
-Euler polynomials with weight zero, $\widetilde{E}_{n,q}\left( x\right) $,
which is defined by Araci $et$ $al.$ and Kim and Choi in \cite{Araci 5} and 
\cite{Kim 7} as follows:%
\begin{equation}
\frac{\widetilde{G}_{n+1,q}\left( x\right) }{n+1}=\lim_{m\rightarrow \infty }%
\frac{1}{\left[ p^{m}\right] _{-q}}\sum_{y=0}^{p^{m}-1}\left( -1\right)
^{y}\left( x+y\right) ^{n}q^{y}  \label{equation 17}
\end{equation}

and%
\begin{equation}
\widetilde{E}_{n,q}\left( x\right) =\int_{%
%TCIMACRO{\U{2124} }%
%BeginExpansion
\mathbb{Z}
%EndExpansion
_{p}}\left( x+y\right) ^{n}d\mu _{-q}\left( y\right) .  \label{equation 18}
\end{equation}

By expressions of (\ref{equation 16}), (\ref{equation 17}) and (\ref%
{equation 18}), we easily discover the following corollary.

\begin{corollary}
For $n\in 
%TCIMACRO{\U{2115} }%
%BeginExpansion
\mathbb{N}
%EndExpansion
^{\ast }$, then we have%
\begin{equation*}
\frac{\left( -1\right) ^{n}}{\left( 1+q\right) ^{n}}\mathcal{A}_{n,\chi
}\left( -q\right) =\frac{d^{n}}{\left( n+1\right) \left[ d\right] _{-q^{-1}}}%
\sum_{a=0}^{d-1}\left( -1\right) ^{a}\chi \left( a\right) q^{-a}\widetilde{G}%
_{n+1,q^{-d}}\left( \frac{a}{d}\right)
\end{equation*}%
Moreover,%
\begin{equation*}
\frac{\left( -1\right) ^{n}}{\left( 1+q\right) ^{n}}\mathcal{A}_{n,\chi
}\left( -q\right) =\frac{d^{n}}{\left[ d\right] _{-q^{-1}}}%
\sum_{a=0}^{d-1}\left( -1\right) ^{a}\chi \left( a\right) q^{-a}\widetilde{E}%
_{n,q^{-d}}\left( \frac{a}{d}\right) \text{.}
\end{equation*}
\end{corollary}

\section{\textbf{On the Eulerian-}$L$\textbf{\ function}}

Our goal in this section is to introduce Eulerian-$L$ function by applying
Mellin transformation to the generating function of Dirichlet's type of
Eulerian polynomials. By (\ref{equation 12}), for $s\in 
%TCIMACRO{\U{2102} }%
%BeginExpansion
\mathbb{C}
%EndExpansion
$, we define the following%
\begin{equation*}
L_{E}\left( s\mid \chi \right) =\frac{1}{\Gamma \left( s\right) }%
\int_{0}^{\infty }t^{s-1}\mathcal{F}_{q}\left( t\mid \chi \right) dt
\end{equation*}

where $\Gamma \left( s\right) $ is the Euler Gamma function. It becomes as
follows:%
\begin{eqnarray*}
L_{E}\left( s\mid \chi \right) &=&q\left[ 2\right] _{q}\sum_{m=0}^{\infty
}\left( -1\right) ^{m}\chi \left( m\right) q^{-m}\left\{ \frac{1}{\Gamma
\left( s\right) }\int_{0}^{\infty }t^{s-1}e^{-m\left( 1+q\right) t}dt\right\}
\\
&=&\frac{q}{\left( 1+q\right) ^{s-1}}\sum_{m=1}^{\infty }\frac{\left(
-1\right) ^{m}\chi \left( m\right) }{q^{m}m^{s}}
\end{eqnarray*}

So, we give definition of Eulerian $L$-function as follows:

\begin{definition}
For $s\in 
%TCIMACRO{\U{2102} }%
%BeginExpansion
\mathbb{C}
%EndExpansion
$, then we have%
\begin{equation}
L_{E}\left( s\mid \chi \right) =\frac{q}{\left( 1+q\right) ^{s-1}}%
\sum_{m=1}^{\infty }\frac{\left( -1\right) ^{m}\chi \left( m\right) }{%
q^{m}m^{s}}\text{.}  \label{equation 19}
\end{equation}
\end{definition}

Substituting $s=-n$ into (\ref{equation 13}), then, relation between
Eulerian $L$-function and Dirichlet's type of Eulerian polynomials are given
by the following theorem.

\begin{theorem}
The following equality holds true:%
\begin{equation*}
L_{E}\left( -n\mid \chi \right) =\left\{ 
\begin{array}{cc}
-\mathcal{A}_{n,\chi }\left( -q\right) & \text{if }n\text{ odd,} \\ 
\mathcal{A}_{n,\chi }\left( -q\right) & \text{if }n\text{ even.}%
\end{array}%
\right.
\end{equation*}
\end{theorem}

%
%%%%%%%%%%%%%%%%%%%%%%%%%%%%%%%%%%%%%%%%%%%%%%%%%%%%%%%%%%%%%%%%

%%%%%%%%%%%%%%%%%%%%%%%%%%%%%%%%%%%%%%%%%%%%%%%%%%%%%%%%%%%%%%%%%%%


\begin{thebibliography}{99}
\bibitem{KIM} D. S. Kim, T. Kim, W. J. Kim and D. V. Dolgy, A note on
Eulerian polynomials, Abstract and Applied Analysis (Article in press).

\bibitem{KIM1} T. Kim, On explicit formulas of $p$-adic $q$-$L$-functions,
Kyushu. J. Math. 48 (1994), 73-86.

\bibitem{KIM2} T. Kim, On a $q$-analogue of the $p$-adic log gamma functions
and related integrals, J. Number Theory 76 (1999), 320-329.

\bibitem{KIM3} T. Kim and S. H. Rim, A note on $p$-adic Carlitz's $q$%
-Bernoulli numbers, Bull. Austral. Math. Soc. Vol. 62 (2000), 227-234.

\bibitem{KIM4} T. Kim, $q$-Volkenborn integration, Russ. J. Math Phys., 19
(2002), 288-299.

\bibitem{KIM5} T. Kim, Non-archimedean $q$-integrals associated with
multiple Changhee $q$-Bernoulli Polynomials, Russ. J. Math Phys., 10 (2003),
91-98.

\bibitem{KIM6} T. Kim, $p$-adic $q$-integrals associated with the
Changhee-Barnes' $q$-Bernoulli Polynomials, Integral Transform. Spec.
Funct., 15 (2004), 415-420.

\bibitem{KIM 6} T. Kim, Analytic continuation of multiple $q$-zeta functions
and their values at negative integers, Russ. J. Math Phys., 11 (2004), 71-76.

\bibitem{Kim 1} T. Kim, $q$-Generalized Euler Numbers and polynomials,
Russian Journal of Mathematical Physics, Vol. 13, No. 3, 2006, pp. 293-308.

\bibitem{Kim 2} T. Kim, Some identities on the $q$-Euler polynomials of
higher order and $q$-stirling numbers by the fermionic $p$-adic integral on $%
%TCIMACRO{\U{2124} }%
%BeginExpansion
\mathbb{Z}
%EndExpansion
_{p}$, Russian J. Math. Phys. 16 (2009), 484--491.

\bibitem{Kim 3} T. Kim, On the $q$-extension of Euler and Genocchi numbers,
J. Math. Anal. Appl. 326 (2007) 1458--1465.

\bibitem{Kim 4} T. Kim, On the analogs of Euler numbers and polynomials
associated with $p$-adic $q$-integral on $%
%TCIMACRO{\U{2124} }%
%BeginExpansion
\mathbb{Z}
%EndExpansion
_{p}$ at $q=1$, J. Math. Anal. Appl. 331 (2007) 779--792.

\bibitem{Kim 5} T. Kim, Some Identities on the integral representation of
the product of several $q$-Bernstein-type polynomials, Abstract and Applied
Analysis, Volume 2011, Article ID 634675, 11 pages.

\bibitem{Kim 6} T. Kim, S. H. Lee, H. H. Han and C. S. Ryoo, On the values
of the weighted $q$-Zeta and $L$-functions, Discrete Dynamics in Nature and
Society, Volume 2011, Article ID 476381, 7 pp.

\bibitem{Kim 7} T. Kim and J. Choi, On the $q$-Euler numbers and polynomials
with weight 0, Abstract and Applied Analysis, vol. 2012, Article ID 795304,
7 pages, 2012.

\bibitem{Jang} L. C. Jang, The $q$-analogue of twisted Lerch type Euler Zeta
functions, Bull. Korean Math. Soc. 47 (2010), No. 6, pp. 1181-1188.

\bibitem{Jang 1} L. C. Jang, V. Kurt, Y. Simsek, and S. H. Rim, $q$-analogue
of the $p$-adic twisted $l$-function, Journal of Concrete and Applicable
Mathematics, vol. 6, no. 2, pp. 169--176, 2008.

\bibitem{Ozden} H. Ozden, I. N. Cangul, Y. Simsek, Multivariate
interpolation functions of higher order $q$-Euler numbers and their
applications, Abstract and Applied Analysis 2008 (2008), Article ID 390857,
16 pages.

\bibitem{Cangul} I. N. Cangul, H. Ozden, Y. Simsek, A new approach to $q$%
-Genocchi numbers and their interpolation functions, Nonlinear Analysis 71
(2009), pp. 793-799.

\bibitem{Cetin} E. Cetin, M. Acikgoz, I. N. Cangul and S. Araci, A note on
the ($h$,$q$)-Zeta type function with weight $\alpha $,
http://arxiv.org/pdf/1206.5299.pdf.

\bibitem{Araci} S. Araci, M. Acikgoz, K. H. Park and H. Jolany, On the
unification of two families of multiple twisted type polynomials by using $p$%
-adic $q$-integral on $%
%TCIMACRO{\U{2124} }%
%BeginExpansion
\mathbb{Z}
%EndExpansion
_{p}$ at $q=-1$, Bulletin of the Malaysian Mathematical Sciences and Society
(Article in press).

\bibitem{Araci 1} S. Araci, M. Acikgoz and K. H. Park, A note on the $q$%
-analogue of Kim's $p$-adic $\log $ gamma type functions associated with $q$%
-extension of Genocchi and Euler numbers with weight $\alpha $, Bulletin of
the Korean Mathematical Society (accepted for publication).

\bibitem{Araci 3} S. Araci, D. Erdal and J. J. Seo, A study on the fermionic 
$p$-adic $q$-integral representation on $%
%TCIMACRO{\U{2124} }%
%BeginExpansion
\mathbb{Z}
%EndExpansion
_{p}$ associated with weighted $q$-Bernstein and $q$-Genocchi polynomials,
Abstract and Applied Analysis, Volume 2011, Article ID 649248, 10 pages.

\bibitem{Araci 4} S. Araci, J. J. Seo and D. Erdal, New construction
weighted $\left( h,q\right) $-Genocchi numbers and polynomials related to
Zeta type function, Discrete Dynamics in Nature and Society, Volume 2011,
Article ID 487490, 7 pages, doi:10.1155/2011/487490.

\bibitem{Araci 5} S. Araci, M. Acikgoz and Feng Qi, On the $q$-Genocchi
numbers and polynomials with weight zero and their applications,
http://arxiv.org/abs/1202.2643.

\bibitem{Simsek} Y. Simsek, $q$-analogue of twisted $l$-series and $q$%
-twisted Euler numbers, Journal of Number Theory, vol. 110, no. 2, pp.
267--278, 2005.

\bibitem{Simsek 1} Y. Simsek, $q$-Dedekind type sums related to $q$-zeta
function and basic $L$-series, Journal of Mathematical Analysis and
Applications, vol. 318, no. 1, pp. 333--351, 2006.

\bibitem{Simsek 2} Y. Simsek, Twisted $\left( h,q\right) $-Bernoulli numbers
and polynomials related to twisted $\left( h,q\right) $-zeta function and $L$%
-function, Journal of Mathematical Analysis and Applications, vol. 324, no.
2, pp. 790--804, 2006.

\bibitem{Acikgoz} M. Acikgoz and Y. Simsek, On Multiple Interpolation
functions oof the N\"{o}rlund-Type $q$-Euler polynomials, Abstract and
Applied Analysis, Volume 2009, Article ID 382574, 14 pages.
\end{thebibliography}
\end{document}